\documentclass{amsart}
\usepackage{amssymb}
\usepackage{amsmath}

\newtheorem{theorem}{Theorem}
\newtheorem{Main theorem}{Main theorem}
\newtheorem{prop}{Proposition}
\newtheorem{lemma}{Lemma}
\newtheorem{corollary}{Corollary}

\newcommand{\LL}{\ensuremath{\mathbb{L}}}

\newcommand{\N}{\ensuremath{\mathbb{N}}}
\newcommand{\Z}{\ensuremath{\mathbb{Z}}}
\newcommand{\Q}{\ensuremath{\mathbb{Q}}}
\newcommand{\R}{\ensuremath{\mathbb{R}}}
\newcommand{\C}{\ensuremath{\mathbb{C}}}
\newcommand{\Pro}{\ensuremath{\mathbb{P}}}
\newcommand{\A}{\ensuremath{\mathbb{A}}}

\begin{document}
\title{Motivic generating series for toric surface singularities}
\author[J. Nicaise]{Johannes Nicaise$^\dag$}
\thanks{$\dag$Research Assistant of the Fund for Scientific
Research --
 Flanders (Belgium)(F.W.O.)}
\address{Department of Mathematics\\
Katholieke Universiteit Leuven\\ Celestijnenlaan 200B\\ B-3001
Leuven\\ Belgium}
\email{johannes.nicaise@wis.kuleuven.ac.be}

\begin{abstract}
Lejeune-Jalabert and Reguera computed the geometric Poincar\'e
series $P_{geom}(T)$ for toric surface singularities. They raise
the question whether this series equals the arithmetic Poincar\'e
series. We prove this equality for a class of toric varieties
including the surfaces,
 and construct a counterexample in the general case.
 We also compute the motivic Igusa Poincar\'e series
$Q_{geom}(T)$ for toric surface singularities, using the change of
variables formula for motivic integrals, thus answering a second
question of Lejeune-Jalabert and Reguera's.
 The series
$Q_{geom}(T)$ contains more information than the geometric series,
since it determines the multiplicity of the singularity. In some
sense, this is the only difference between $Q_{geom}(T)$ and
$P_{geom}(T)$.
\end{abstract}

\maketitle
\section{Introduction}
Throughout this article, we work over a base field $k$ of
characteristic zero, and we denote by $k'
$ its algebraic closure.

\noindent The classical motivation for the introduction of the
motivic generating series can be found in $p$-adic analysis. Let
$L$ be a finite field extension of $\Q_p$, with ring of integers
$\mathcal{O}_L$ and uniformizing parameter $\pi$. Let $X$ be a
variety over $\mathcal{O}_L$.

 The Igusa Poincar\'e series counts approximate solutions modulo
 $\pi^{n+1}$. To be precise, let $\tilde{N}_n$ be the number of
 points $|X(\mathcal{O}_L/\pi^{n+1})|$, for $n\geq 0$. Then the
 Igusa Poincar\'e series is defined to be
 \[Q(T)=\sum_{n\geq 0} \tilde{N}_{n} T^{n}\ .\]

 The Serre-Oesterl\'e series counts approximate solutions that can
 be lifted to global solutions on $X$: putting $\bar{N}_{n}$ equal
 to the cardinality of the image of $X(\mathcal{O}_{L})$ in
 $X(\mathcal{O}_L/\pi^{n+1})$, the series is defined as
 \[P(T)=\sum_{n\geq 0}\bar{N}_{n} T^{n}\ .\]

Both series are known to be rational: Igusa proved the rationality
of $Q(T)$ in the hypersurface case, rewriting the series as a
$p$-adic integral and applying resolution of singularities
\cite{Igusa2}. Denef proved the rationality of $P(T)$ making use
of the model-theoretic framework of quantifier elimination and
cell decomposition \cite{Denef}.

Since motivic integration is introduced as a formal analogue of
$p$-adic integration, replacing the ring $\Z_p$ by $k[[t]]$ and
taking values in the completed localized Grothendieck ring
$\hat{\mathcal{M}}_k$, it is natural to translate these series to
the motivic setting. This was done by Denef and Loeser (see e.g.
\cite{DL2}\cite{DL3}). Let $X$ be a variety over $k$. The Igusa
Poincar\'e series has a straightforward motivic counterpart. The
approximate solutions are given by $n$-jets, i.e. points on the
scheme $\mathcal{L}_{n}(X)$, which will be introduced in the next
section. Instead of counting points, we use the universal additive
invariant, mapping a constructible set to its isomorphism class in
the Grothendieck ring. In this way, we obtain
\[Q_{geom}(T)=\sum_{n\geq 0}[\mathcal{L}_{n}(X)]T^{n}\ .\] This
series is rational in $\mathcal{M}_{k}[[T]]$, as can be proven by
making use of resolution of singularities, and the change of
variables formula for motivic integrals. Furthermore, when $X$ is
defined over some number field $L$, this series specializes to the
classical Igusa Poincar\'e series for almost all finite places
$\mathcal{P}$. By this we mean the following: we can choose a
model over $\mathcal{O}_{L}$ for $X$, and count points modulo
$\mathcal{P}$, for each finite place $\mathcal{P}$. This operation
is denoted by an operator $N_{\mathcal{P}}$; $N_{\mathcal{P}}(X)$
is well-defined for almost all finite places $\mathcal{P}$.
Applying $N_{\mathcal{P}}$ termwise to the series $Q_{geom}$
yields the Igusa Poincar\'e series $Q(T)$ for
$X\times$\,Spec\,$\mathcal{O}_{L_{\mathcal{P}}}$, for almost all
finite places $\mathcal{P}$.

A naive generalization for $P(t)$ is obtained by looking at
$n$-jets that can be lifted to arcs on $X$, that is, by defining
the geometric Poincar\'e series as
\[P_{geom}=\sum_{n\geq 0}[j^{n}(\mathcal{L}(X))]T^{n}\ .\] This
series is well-defined, since a theorem of Greenberg guarantees
that $j^{n}(\mathcal{L}(X))$ is constructible, and it is rational
in $\mathcal{M}_{k}[[T]]$ (see \cite{DLinvent}). But, in general,
this series does not specialize to the Serre-Oesterl\'e series
when $X$ is defined over a number field. The reason for this is
that, working scheme-theoretically, we allow extensions of the
base field when lifting jets. So instead of counting approximate
solutions which can be lifted to a solution in
$\mathcal{O}_{L_{\mathcal{P}}}$, we count approximate solutions
that can be lifted to a solution in a maximal unramified extension
of $\mathcal{O}_{L_{\mathcal{P}}}$, whose residue field is
precisely the separable closure of the residue field of
$\mathcal{O}_{L_\mathcal{P}}$. The arithmetic Poincar\'e series is
designed to remedy this discrepancy. While the Igusa Poincar\'e
series can be computed from a resolution of singularities, the
geometric and arithmetic series are very hard to compute. The only
known cases so far are formal branches of plane curves \cite{DL},
toric surfaces \cite{LejReg}, and surfaces with an embedded
resolution of a simple form \cite{Nic2} (this latter result yields
an easy way to recover the formula for toric surfaces).

Let us give an overview of the results in this paper. In Section
\ref{Motivic}, we recall the definitions of the motivic generating
series (Igusa Poincar\'e series, geometric Poincar\'e series,
arithmetic Poincar\'e series), with emphasis on the last one. In
Section \ref{toric}, we develop a sufficient condition for the
equality of the geometric and arithmetic Poincar\'e series of
toric varieties (Theorem \ref{equal}), and we show that this
condition is always satisfied in the case of a toric surface
(Corollary \ref{surf}). Section \ref{counter} contains an example
of a toric threefold for which the series differ (Proposition
\ref{example}). Section \ref{integration} gives a concise
introduction to the theory of motivic integration. In order to
compute the Igusa Poincar\'e series $Q_{geom}(T)$ of a toric surface, we
establish in Section \ref{igusa} a factorization of the minimal
toric resolution into a sequence of blow-ups of smooth subschemes.
The actual computation is done in Section \ref{comput}. In Section
\ref{poles}, we determine which information is contained in the
Igusa Poincar\'e series, by investigating its poles, and we
compare it with the formula for the geometric series $P_{geom}$ in
\cite{LejReg}. Theorem \ref{info} shows that $Q_{geom}(T)$ contains more
information than $P_{geom}$, since $Q_{geom}(T)$ also determines the
multiplicity of the toric surface singularity. In some sense, this
is the only additional information you obtain from $Q_{geom}(T)$.

\section{Motivic Poincar\'e series}\label{Motivic}
Let $X$ be a variety over $k$, that is, a reduced and separated
scheme of finite type over $k$, not necessarily irreducible. For
each positive integer $n$, the functor from the category of
$k$-algebras to the category of sets, sending an algebra $R$ to
the set of $R[[t]]/t^{n+1}R[[t]]$-rational points on $X$, is
representable by a scheme $\mathcal{L}_{n}(X)$. Since the natural
projections $j^{n}_{n+1}:\mathcal{L}_{n+1}(X)\rightarrow
\mathcal{L}_{n}(X)$ are affine, we can take the projective limit
in the category of schemes to obtain the scheme of arcs
$\mathcal{L}(X)$. This scheme represents the functor sending a
$k$-algebra $R$ to the set of $R[[t]]$-rational points on $X$, and
comes with natural projections $j^{n}:\mathcal{L}(X)\rightarrow
\mathcal{L}_{n}(X)$, mapping an arc to its $n$-truncation.
We consider $\mathcal{L}(X)$ and $\mathcal{L}_{n}(X)$ as endowed 
with their reduced structure. By an
arc (resp. $n$-jet) on $X$, we always mean a $k'$-rational
point on $\mathcal{L}(X)$ (resp. on $\mathcal{L}_n(X)$), unless
explicitly stated otherwise. When $X$ is smooth, the morphisms
$j^{n}_{n+1}$ are Zariski-locally trivial fibrations with fiber
$\A^{d}_{k}$, where $d$ is the dimension of $X$.

 We now introduce the Grothendieck ring $K_{0}(Var_k)$ of varieties over $k$.
 Start from the free abelian group generated by isomorphism
 classes $[X]$ of varieties $X$ over $k$, and consider the
 quotient by the relations $[X]=[X\setminus X']+[X']$, where $X'$
 is closed in $X$. A constructible subset of $X$ can be written as a disjoint union of
 locally closed subsets and determines unambiguously an element of $K_{0}(Var_k)$.
  The Cartesian product induces a product on
 $K_{0}(Var_k)$, which makes it a ring. We denote the class of the affine line
 $\A^{1}_{k}$ in $K_{0}(Var_k)$ by $\LL$, and the localization of $K_{0}(Var_k)$ with
 respect to $\LL$ by $\mathcal{M}_k$. On $\mathcal{M}_{k}$, we consider a decreasing filtration $F^{m}$, where $F^{m}$
 is the subgroup generated by elements of the form $[X]\LL^{-i}$, with dim\,$X-i\leq -m$. We define
 $\hat{\mathcal{M}}_{k}$ to be the completion of $\mathcal{M}_{k}$ with respect to this filtration.

 The Grothendieck ring $K_{0}(Var_k)$ is
 not very well understood. Recently, Poonen showed that it is not
 a domain \cite{Poo}.
Bittner proved in \cite{Bitt}, using the Weak Factorization
Theorem, that $K_{0}(Var_k)$
  can be presented by taking the
 isomorphism classes of smooth projective varieties as generators,
 and by considering the relations $[\emptyset]=0$ and
 $[Bl_{Y}X]-[E]=[X]-[Y]$, where $X$ and $Y$ are smooth projective
 varieties, $Y\subset X$, and $E$ is the exceptional divisor of
 the blow-up $Bl_{Y}(X)$ of $X$ along $Y$. This presentation is
 important for the construction of additive invariants; it allows
 one to prove the existence of a ring morphism $\chi_{mot}$ from
 the Grothendieck ring of varieties over $k$ to the Grothendieck
 ring of Chow motives over $k$, sending the class of a smooth
 projective variety to the class of its associated Chow motive,
 and sending $\LL$ to the class of the Tate motive $\LL_{mot}$.
 The existence of this map was proven already in \cite{GuiNav}. We
 denote the image of this morphism by $K_{0}^{mot}(Var_k)$.
 A definition of Chow motives can be found in \cite{Scholl}; the idea is
 that motives should provide some kind of universal cohomology
 theory. Smooth projective varieties with isomorphic Chow motives
 have the same cohomology for all known cohomology theories with
 coefficients in a field of characteristic zero.

The motivic Igusa Poincar\'e series is defined as
\[Q_{geom}(T)=\sum_{n\geq 0}[\mathcal{L}_{n}(X)]T^{n},\] while the
geometric Poincar\'e series is by definition
\[P_{geom}(T)=\sum_{n\geq 0}[j^{n}(\mathcal{L}(X))]T^{n}\ .\] The
latter series is well defined, since Greenberg's theorem \cite{Gr}
states that we can find a positive integer $c$ such that, for all
$n$, and for each field $K$ containing $k$,
$j^{n}(\mathcal{L}(X)(K))=j^{n}_{nc}(\mathcal{L}_{nc}(X)(K))$. So
it follows from Chevalley's theorem \cite{Hart} that
$j^{n}(\mathcal{L}(X))$ is constructible, and hence determines an
element $[j^{n}(\mathcal{L}(X))]$ in $K_{0}(Var_k)$. One can
define local variants of both series by only considering arcs with
origin in some closed subvariety $Z$ of $X$, as is done in
\cite{LejReg}. If we write $\mathcal{L}(X)_{Z}$ to denote the
inverse image $(j^0)^{-1}(Z)$ in $\mathcal{L}(X)$, the local
geometric Poincar\'e series of $X$ at $Z$ is defined as
\[P_{geom,Z}(T)=\sum_{n\geq
0}[j^{n}(\mathcal{L}(X)_{Z})]T^{n}\,,\] and the local Igusa
Poincar\'e series is defined analogously.

 As mentioned before, the main characteristic of $P_{arith}$ should
 be that it behaves well under specialization to
 $\mathcal{P}$-adic completions of a number field $L$. The crucial
 point in the construction is the use of pseudo-finite fields. A
 pseudo-finite field is an infinite perfect field with exactly one
 field extension of any given finite degree, and over which every
 absolutely irreducible variety has a rational point. Their
 relevance is illustrated by the following theorem of Ax \cite{Ax}: two ring
 formulas over $\Q$ are equivalent when interpreted in
 $\mathbb{F}_p$, for all sufficiently large primes $p$, if and
 only if they are equivalent when interpreted in $K$, for all
 pseudo-finite fields $K$ containing $\Q$. In this way, they
 present themselves as natural candidates to control the behaviour
 of $P_{arith}$ under specialization. Intuitively, they are
 perfectly suited to detect rationality conditions of the form
 $(\exists y)y^{n}=x$, since, for $n>1$ and $p$ sufficiently large, not every element of $\mathbb{F}_p$ can have an $n$-th root.
 This means that the condition $(\exists y)y^{n}=x$, which is
 ignored when working over an algebraically closed field, is
 always brought into account by the much more sensitive
 pseudo-finite fields. The counterexample in section 4, and the
 exact definition in \cite{DL1} of the map $\chi_c$ introduced
 below, will clarify this remark.

  A ring formula over a
 field $k$ is a logical formula $\varphi$ built from Boolean combinations
 of polynomial equalities over $k$, and quantifiers. When allowing
 extension to an algebraically closed field, we can eliminate quantifiers from
 $\varphi$ and thus associate to $\varphi$ an element of the
 Grothendieck ring. But these field extensions are exactly what we
 try to avoid. We will associate to $\varphi$ an element of
 $K_{0}^{mot}(Var_k)\otimes \Q$ in a more subtle way.

Consider the Grothendieck ring $K_{0}(PFF_k)$ of the theory of
pseudo-finite fields containing $k$. It is generated by classes
$[\varphi]$, where $\varphi$ is a ring formula over $k$, which are
subject to the relations
$[\varphi_1\vee\varphi_2]=[\varphi_1]+[\varphi_2]-[\varphi_1
\wedge \varphi_2]$, whenever $\varphi_1$ and $\varphi_2$ have the
same free variables, and to the relations
$[\varphi_1]=[\varphi_2]$, whenever there exists a ring formula
$\psi$ over $k$ such that, interpreted over any pseudo-finite
field $K$ containing $k$, $\psi$ defines a bijection between the
tuples over $K$ satisfying $\varphi_1$ and those satisfying
$\varphi_2$. Ring multiplication is induced by taking the
conjunction of formulas in disjoint sets of variables. Denef and
Loeser \cite{DL1} constructed a morphism
\[\chi_{c}:K_{0}(PFF_k)\rightarrow K_{0}^{mot}(Var_k)\otimes \Q\
.\] For this construction, it is important to understand the
structure of $K_{0}(PFF_k)$. The theory of quantifier elimination
for pseudo-finite fields \cite{FriJar}\cite{FriSac}, states that
quantifiers can be eliminated if one adds some relations to the
language, which have a geometric interpretation in terms of Galois
covers. This interpretation yields a construction for $\chi_{c}$.
It is important for our purposes that, if our original ring
formula $\varphi$ did not contain any quantifiers in the first
place, $\chi_{c}$ maps $[\varphi]$ to the class of the
constructible set defined by $\varphi$ in $K^{mot}_{0}(Var_k)$.

Now, we are ready to define the arithmetic Poincar\'e series
$P_{arith}$. We only consider the case where $X$ is a subvariety
of some affine space $\A^{m}_k$; the general case can be dealt
with using definable subassignements \cite{DL1}. It follows from
Greenberg's theorem that we can find, for each positive integer
$n$, a ring formula $\varphi_{n}$ over $k$, such that, for all
fields $K$ containing $k$, the $K$-rational points of
$\mathcal{L}_{n}(X)$ that can be lifted to a $K$-rational point of
$\mathcal{L}(X)$, correspond to the tuples satisfying the
interpretation of $\varphi_{n}$ in $K$. We define the arithmetic
Poincar\'e series to be
\[P_{arith}=\sum_{n\geq 0}\chi_{c}([\varphi_{n}])T^{n}\ .\]
The local variant $P_{arith,Z}$, where $Z$ is a closed subvariety
of $X$ (defined over $k$), is defined in the obvious way.

 The series is
rational over $K_{0}^{mot}(Var_k)[\LL^{-1}]\otimes\Q$. If $k$ is a
number field $L$ , we recover the Serre-Oesterl\'e series for
$X\times$\,Spec\,$\mathcal{O}_{L_{\mathcal{P}}}$, for almost all
finite places $\mathcal{P}$, by applying the operator
$N_{\mathcal{P}}$ to each coefficient of numerator and denominator
\cite{DL1}.

\section{The arithmetic Poincar\'e series for toric
varieties}\label{toric} Proposition 3.3 in \cite{LejReg}
identifies an arc through the zero-dimensional orbit $O$ on an
affine toric surface $X$ meeting the embedded torus $T$ with a
couple consisting of an arc on $T$ and an $N$-vector in the
interior of the cone $\sigma$ associated to $X$. We sketch this
identification, which generalizes immediately to arbitrary
dimensions. Let $X$ be an affine toric variety, defined over a
field $k$ with algebraic closure $k'$, and associated to an
$n$-dimensional cone $\sigma$ in $N_{\R}=N\otimes\R$, where $N$ is
a lattice of dimension $n$. Let $O$ be the unique orbit of
dimension zero, and let $T$ be the orbit of dimension $n$. We
denote by $M$ the dual lattice of $N$, and by $\check{\sigma}$ the
dual cone of $\sigma$ in $M_{\R}=M\otimes \R$.

Let $h$ be an arc on $X$ through $O$ meeting $T$. This arc can be
represented by a coordinate morphism $\psi_{h}:\check{\sigma}\cap
M\rightarrow k'[[t]]$. Since $h$ meets $T$, the image of this
morphism does not contain zero, and thus we can define a new
mapping $\check{\sigma}\cap M\rightarrow \N$ by composing with the
function $ord_{t}$, measuring the order of a power series in
$k'[[t]]$. This mapping extends to a linear form
$\nu_{h}:M\rightarrow \Z$, defining a vector in $N$, which is
contained in the interior $Int(\sigma)$ of $\sigma$, since
$\nu_{h}(m)>0$ whenever $m\in \check{\sigma} \cap M$ and $m\neq
0$. If we set $u_{h}(m)=\psi_{h}(m)t^{-\nu_h(m)}$ for $m\in
\check{\sigma}\cap M$, the mapping $u_{h}$ extends to a morphism
from $M$ to the multiplicative group of units in $k'[[t]]$, which
is nothing but an arc on $T$. Conversely, $h$ can be recovered
from $\nu_{h}$ and $u_{h}$ by setting
$\psi_{h}(m)=t^{<m,\nu_{h}>}u_{h}(m)$.

One could say that we have split up the arc $h$ into an order
function $\nu_{h}$ and an angular component $u_{h}$. In what
follows, we identify an arc $h$ with the associated couple
$(\nu_h,u_h)$. The smoothness of $T$ reduces the computation of
the geometric Poincar\'e series to a combinatorial analysis of the
behaviour of $\nu_{h}$ when $h$ varies.

The next thing we have to do, is to prove that the arcs meeting
$T$ suffice to compute the motivic Poincar\'e series. Let $H$ be
the set of arcs through $O$ on $X$, and let $H^{*}$ be the subset
consisting of arcs meeting $T$.

\begin{lemma}[Moving Lemma]
$j^{s}(H)=j^{s}(H^{*})$ for each $s$.
\end{lemma}
\begin{proof}
 Let $h:\tilde{X}\rightarrow X$ be a toric resolution of $X$, corresponding to a subdivision of $\sigma$ into a simple fan $\Sigma$.
   Let $\psi$ be an arc on $X$. We will prove
that we can deform $\psi$ to an arc $\psi'$ meeting $T$ without
changing its $s$-jet. Let $\eta$ be the image of the generic point
of Spec\,$k'[[t]]$ under $\psi$, and let $\tau$ be the face of
$\sigma$ such that $\eta$ is contained in the orbit $O_{\tau}$
corresponding to $\tau$. Since the choice of a cone of $\Sigma$,
contained in $\tau$ and of the same dimension, yields a section of
$h$ over $O_{\tau}$, we can lift
 the morphism
Spec\,$k'((t))\rightarrow X$ induced by $\psi$ to a morphism
Spec\,$k'((t))\rightarrow \tilde{X}$.

Applying the valuative criterion for properness to the morphism
$h$, we see that the morphism Spec\,$k'((t))\rightarrow \tilde{X}$
 has a unique extension to
an arc $\tilde{\psi}$ on $\tilde{X}$. It is clear that
$h(\tilde{\psi})=\psi$. The variety $\tilde{X}$ being smooth, it
is easy to move $\tilde{\psi}$ away from $h^{-1}(X-T)$ (i.e. out
of the inverse image of $\mathcal{L}(X-T)$ under $h$) without
changing its $s$-jet, using a system of local parameters. Now take
$\psi'$ to be the image under $h$ of the arc $\tilde{\psi}'$
obtained in this way.
\end{proof}

 To prove the equality of the series $P_{arith}$ and $P_{geom}$, we have to find a convenient
way to describe truncations of an arc $h=(\nu_{h},u_{h})$ in
$H^{*}$. Let $\tau$ be a face of $\sigma$, denote by $N_{\tau}$
the sublattice of $N$ generated by $\tau\cap N$, and let
$M_{\tau}$ be its dual. Let $\check{G}_{\tau}$ be a minimal set of
generators for the semigroup $\check{\tau}\cap M_{\tau}$. Suppose
that we can find, for each face $\tau$ of $\sigma$, and for each
vector $\nu_h$ in $N\cap Int(\tau)$, a basis
$\{\mu_{i}\}_{i=1}^{\mathrm{dim}\,\tau}$ for $M_{\tau}$,
consisting of elements of $\check{G}_{\tau}$, such that
$<\mu,\nu_{h}>\,\geq\, <\mu_{i},\nu_{h}>$ for each $\mu$ in
$\check{G}_{\tau}\setminus
\{\mu_{i}\}_{i=1}^{\mathrm{dim}\,\tau}$, and for each $i$ such
that $\mu$ is not contained in the coordinate hyperplane
$\lambda_{i}=0$ defined by $\mu_{i}$ (*).

\begin{theorem}\label{equal}
If (*) holds, then $P_{geom}=P_{arith}$ in
$K_0^{mot}(Var_k)\otimes \Q$.
\end{theorem}
\begin{proof}
 Because of the torus action on $X$, the global series
$P_{geom}$ and $P_{arith}$ can be written in terms of the local
series $P_{geom,x_{\tau}}$ and $P_{arith,x_{\tau}}$ at the
distinguished point $x_{\tau}$ of $O_{\tau}$, where $\tau$ is a
face of $\sigma$. To be precise,
\[P_{geom}=\sum_{\tau\leq \sigma}(\LL-1)^{n-\mathrm{dim}\,
\tau}P_{geom,x_{\tau}}\,,\] and the analogous statement holds for
$P_{arith}$. Hence, it suffices to prove the theorem for the local
series at $x_{\tau}$. If we denote by $Y$ the complement in $X$ of
all orbits $O_{\tau'}$ with $\tau'$ a face of $\sigma$ that is not
contained in $\tau$, then $Y$ is isomorphic to the toric variety
associated to the cone $\tau$ in $N_{\R}$. Let $N'$ be the
sublattice of $N$ generated by $\tau\cap N$, and let $Y'$ be the
toric variety associated to the cone $\tau$ in $N'_{\R}$. Since
$Y$ is isomorphic to the product of $Y'$ with a torus, it suffices
to prove the theorem for the local series at the zero-dimensional
orbit $O$ of $X$.

 Let $h=(\nu_h,u_h)$ be an arc in $H^{*}$. The angular component $u_{h}$ is completely determined by
$u_{h}(\mu_{i})$, $i=1,\ldots,n$, where $\{\mu_i\}_{i=1}^{n}$ is a
basis of $M$, satisfying (*) for $\tau=\sigma$ and for the order
vector $\nu_h$. Let $s$ be a positive integer. We define a new
angular component $u'_{h}$, mapping $\mu_{i}$ to the truncation of
$u_{h}(\mu_{i})$ at $t^{s+1-<\mu_{i},\nu_{h}>}$ if
$s\,\geq\,<\mu_{i},\nu_{h}>$, and to zero in the other case.
Because of our supposition, $j^{s}(\nu_{h},u_{h})$ equals
$j^{s}(\nu_{h},u'_{h})$. Moreover, $(\nu_{h},u'_{h})$ has
coefficients in the same field as its $s$-jet, since $u'_{h}$ can
be recovered from the $\mu_{i}$-coordinates of this jet. So the
ring formula $\varphi_{s}$ over $k$, used to define $P_{arith}$,
whose interpretation in any field $K$ containing $k$ is the
condition of liftability of a $K$-rational point of
$\mathcal{L}_{s}(X)$ to a $K$-rational point of $\mathcal{L}(X)$,
is equivalent to the set of equalities and inequalities describing
the constructible set $j^{s}(H)$. By definition of the morphism
$\chi_{c}:K_{0}(PFF_k)\rightarrow K_0^{mot}(Var_k)\otimes \Q$,
this implies that the image of $[j^s(H)]$ in
$K_0^{mot}(Var_k)\otimes \Q$ equals $\chi_c(\varphi_s)$, which
proves the equality of the geometric and the arithmetic Poincar\'e
series.
\end{proof}
An explicit expression for $P_{geom}$ is known only in the case
$n=2$. It follows from the proof that, in order to compute the
geometric or arithmetic series for toric varieties of higher
dimension, it suffices to consider the local series at the
zero-dimensional orbit.

\begin{corollary}\label{surf}
The geometric and arithmetic Poincar\'e series of a toric surface
singularity coincide.
\end{corollary}

\begin{proof}
Let $\nu_h$ be a vector in $N\cap Int(\sigma)$, and let
$\check{G}$ be the minimal set of generators of the semi-group
$\check{\sigma}\cap M$. It is easy to see that there exist
elements $\mu_1,\mu_2$ in $\check{G}$,
 forming a $\Z$-basis of $M$, such that $<\mu,\nu_h>\,\geq\, <\mu_{i},\nu_{h}>$ for each $\mu$ in
$\check{G}\setminus\{\mu_1,\mu_2\}$ and $i=1,2$; thus condition
(*) is always satisfied in the surface case.
\end{proof}

This basis is used in \cite{LejReg} to calculate the image under
the truncation map $j^{s}$ of the set $H^{*}_{\nu}$, consisting of
all arcs in $H^{*}$ with fixed order vector $\nu$.

\section{A counterexample}\label{counter}
Our supposition (*) is valid when $n=2$, but does not always hold
in higher dimensions.

\begin{prop}\label{example}
There exists an affine toric threefold $V$, with zero-dimensional
orbit $O$, such that the local geometric and arithmetic Poincar\'e
series at $O$ differ.
\end{prop}

\begin{proof}
Consider the cone $\check{\sigma}$ generated by $(1,0,0),
(0,1,0)$, and $(1,1,2)$, and put $\nu$ equal to $(2,2,-1)$. There
are three lattice points in $\check{\sigma}$, minimizing $\nu$,
and these are precisely the generators - which do not form a basis
for $M$.

The problem that arises is the following: our angular component
$u_h$ is determined by its values at the basis
$\{(1,0,0),\,(0,1,0),\,(1,1,1)\}$. To compute the $(1,1,1)$-
coordinate of the $s$-jet of $u_h$, we only need the first $s-3$
coefficients of $u_h(1,1,1)$; but in order to compute the
$(1,1,2)$-coordinate, we also need the $(s-2)$-th. So it may
happen that these first $s-3$ coefficients lie in a field $k$,
while the only $(s-2)$-th coefficient yielding the right value for
the jet of $u_h(1,1,2)$ lies in $k'\setminus k$.

Consider, for instance, the 2-jet mapping $(1,0,0)$ and $(0,1,0)$
to $t^{2}$, $(1,1,1)$ to $0$, and $(1,1,2)$ to $-t^{2}$. If we
define the angular component $u$ by $u(1,0,0)=u(0,1,0)=1$ and
$u(1,1,1)=i$, then $(u,\nu)$ lifts this jet to an arc over $\C$.
However, the jet is not liftable to a $\Q$-rational point of the
arc space, since such an arc $(u',\nu)$ has to satisfy
$u'(1,1,1;0)^{2}=u'(1,1,2;0)u'(1,0,0;0)u'(0,1,0;0)=-1$.

Of course, this does not necessarily mean that this discrepancy
actually emerges in the $T^{2}$-coefficients of the series; so let
us make explicit computations. The set
$\{\mu_1=(1,0,0),\mu_2=(0,1,0),\mu_3=(1,1,1),\mu=(1,1,2)\}$ is a
set of generators for the semigroup $\check{\sigma}\cap M$, and
the first three of them form a lattice basis. An arc $h$, meeting
the embedded torus, is, as always, determined by a vector $\nu$ in
$Int(\sigma)\cap N$, and an angular component $u$. What will its
2-jet look like? For each $i=1\ldots 3$, $u(\mu_i;j)$ can take
random values for $0\leq j\leq 2-<\mu_i,\nu>$. Because (*) is not
satisfied, even when we fix these values, we will have some
freedom in the choice of the coefficients $u(\mu;j)$, for $0\leq
j\leq 2-<\mu,\nu>$.

 If $<\mu,\nu>\,>2$, we
have no worries, and if $<\mu,\nu>$ is greater than or equal to
the maximum of the $<\mu_i,\nu>$, our 2-jet is fixed by the
choices we made. So let us suppose it is not. We may as well
assume that $<\mu_3,\nu>$ is strictly maximal among the
$<\mu_i,\nu>$, because otherwise, $\nu$ satisfies (*). If
$<\mu,\nu>=2$, $u_h(\mu,0)$ is arbitrary, if we allow liftings
over an algebraic closure of $k$. The equalities $<\mu,\nu>=0$ and
$<\mu,\nu>=1$ cannot occur under the assumptions we made.

 This means that the class of $j^{2}(\mathcal{L}(V)_{O})$ in the
 Grothendieck ring $K_{0}(Var_k)$, where $V$
 is the affine toric threefold with zerodimensional orbit $O$
 defined by the dual cone $\sigma$ of $\check{\sigma}$, and where $\mathcal{L}(V)_{O}$ denotes the space of
 arcs on $V$ with origin at $O$, is equal to
 \[\LL^{9}-\LL^{6}+3\LL^{5}-6\LL^{4}+10\LL^{3}-9\LL^{2}+3\LL\,.\]

 \noindent As for the arithmetic series, the only jets liftable over the
 algebraic closure of $k$, but not necessarily over $k$ itself, are the
 2-truncations of arcs $(\nu,u)$ with
 \[<\mu,\nu>=<\mu_1,\nu>=<\mu_2,\nu>=2\,,\] as is the case in our
 example above, and with $u(\mu_1;0)u(\mu_2;0)u(\mu;0)$ not a
 square in $k$. So our ring formula $\varphi_{2}$ becomes
 \[\psi\wedge (\exists y)y^{2}=x_{\mu_1,2}\,x_{\mu_2,2}\,x_{\mu,2}\,,\]
 where $\psi$ is the quantifier-free ring formula over $k$ describing
 $j^{2}(\mathcal{L}(V)_{O})$, and \[x_{\mu_i,j},\,x_{\mu,j},
 \quad i=1,2,3,\,j=0,1,2\] are coordinates of 2-jets in the ambient space
 $\A^{4}$ associated to our set of generators. Rewriting the
 formula as a strict disjunction, we see that, in order to prove
 that the $T^{2}$-coefficients in $P_{geom}$ and $P_{arith}$
 differ, we have to prove the inequality
$\chi_{c}(\varphi'_{2})\neq (\LL-1)^{3}$ in
$K_{0}^{mot}(Var_k)\otimes \Q$,
 where $\varphi'_2$ is the ring formula expressing
\begin{eqnarray*}&&(\forall i)(\forall j\neq
2)(x_{\mu_i,j}=0\,\wedge\, x_{\mu,j}=0)\,\wedge\,
x_{\mu_3,2}=0\,\wedge\, x_{\mu_1,2}\neq
 0 \\&& \qquad \qquad \qquad \qquad\wedge\,
x_{\mu_2,2}\neq 0 \,\wedge\, x_{\mu,2}\neq 0 \,\wedge\, (\exists
y)y^{2}=x_{\mu_1,2}\,x_{\mu_2,2}\,x_{\mu,2}\,,
 \end{eqnarray*}
and where we abuse notation by writing $\LL$ for the class of the
Tate motive in $K_{0}^{mot}(Var_k)\otimes \Q$.

 Let $T$ be the
threefold torus $(\A_{k}^{1}\setminus 0)^{3}$, and consider the
Galois cover
\[f:\mathrm{Spec}\,k[t_i,t_i^{-1},w]/(w^{2}-t_1t_2t_3)\rightarrow T:(t_1,t_2,t_3,w)\mapsto
(t_1,t_2,t_3)\] with Galois group $\Z_2$. For each field $K$
containing $k$, the $K$-rational points of $T$ that lift to a
$K$-rational point of $T$ with respect to $f$, are exactly the
tuples $(t_1,t_2,t_3)$ in $K^{3}$ satisfying $(\exists
y)(y^{2}=t_1t_2t_3)$. By definition of the morphism $\chi_{c}$,
this implies that
\[\chi_{c}(\varphi'_2)=\frac{1}{2}(\LL-1)^{3}\,,\] which is, of course, exactly what one would expect. Since a threefold torus has
non-trivial cohomology, $1/2(\LL-1)^{3}$ cannot be zero in
$K_{0}^{mot}(Var_k)\otimes\Q$, so $P_{arith}\neq P_{geom}$.
\end{proof}

\section{Motivic integration}\label{integration}
 This section gives a concise survey of some definitions concerning
 motivic integration. More exact statements and proofs can be
 found in \cite{DLinvent}.

Motivic integration is a wonderful theory, which may be considered
as an analogue of $p$-adic integration, replacing the ring $\Z_p$
by $k[[t]]$ and taking values in the completed localized
Grothendieck ring $\hat{\mathcal{M}}_k$. Batyrev \cite{Baty} used
$p$-adic integration and the Weil conjectures to prove that
birationally equivalent Calabi-Yau varieties have the same Betti
numbers. Kontsevich \cite{Ko} observed that the development of a
geometrical analogue of $p$-adic integration would allow one to
prove stronger results:
 he used motivic integration to prove that birationally equivalent
Calabi-Yau varieties have the same Hodge numbers. Motivic
integration was further developed by Batyrev
\cite{Baty1}\cite{Baty2}, and Denef and Loeser
\cite{DLinvent}\cite{DL6}\cite{DL4}.

Let us first introduce a motivic measure $\mu$ on the arc space of
a variety $X$ of pure dimension $d$ over $k$. Let $A$ be a subset
of $\mathcal{L}(X)$. We call $A$ a cylinder if there exists a
positive integer $n$, and a constructible subset $A_{n}$ of
$\mathcal{L}_{n}(X)$, such that $A=(j^{n})^{-1}(A_{n})$. We say
$A$ is stable at level $n$ if furthermore, for each $m\geq n$, the
projection $j^{m}_{m+1}$ is a locally trivial fibration over
$j^{m}(A)$, with fiber $\A^{d}$. In this case, we define
$\tilde{\mu}(A)$ to be $[A_{n}]\LL^{-(n+1)d}$ (some authors use
$[A_{n}]\LL^{-nd}$ instead). If $A$ is a cylinder, not necessarily
stable, we define $\mu(A)$ by cutting out tubular neighbourhoods
of $\mathcal{L}(X_{sing})$ in order to obtain stability:
\[\mu(A)=\lim_{e\to \infty}\tilde{\mu}\{A
\setminus(j^{e})^{-1}j^{e}(\mathcal{L}(X_{sing}))\}\ .\] One can
check that this limit exists in $\hat{\mathcal{M}}_{k}$. All these
definitions are inspired by the $p$-adic case. The measure $\mu$
is $\Sigma$-additive: if a cylinder $A$ can be written as a union
of cylinders $A_{i}$, $i\in \N$, then $\mu(A)=\sum_{i}\mu(A_{i})$.

We can extend our class of measurable subsets in the following
way: consider the norm function
$\|.\|:\hat{\mathcal{M}}_{k}\rightarrow \R_{\geq 0}$, mapping an
element $x$ to $2^{-n}$, where $x \in F^{n}$ and $x \notin
F^{n+1}$. We call a subset $A$ of $\mathcal{L}(X)$ measurable if
we can find, for each $\epsilon>0$, a collection $A_{i}(\epsilon)$
of cylinders, $i\in \N$, such that the symmetric difference of $A$
and $A_{0}(\epsilon)$ is contained in $\cup_{i\geq
1}A_{i}(\epsilon)$, and $\|\mu(A_{i}(\epsilon)\|\leq \epsilon$ for
each $i\geq 1$. One can prove that in this case,
$\mu(A)=\lim_{\epsilon\to 0}\mu(A_{0}(\epsilon))$ is well-defined.
So we define the measure of a measurable set by approximating it,
using cylinders.

Let $A\subset \mathcal{L}(X)$ be a measurable set, and let
$\alpha$ a function from $A$ to $\Z\cup\{\infty\}$. We say
$\LL^{-\alpha}$ is integrable if $\alpha^{-1}(i)$ is measurable,
for each $i\in \Z$, and if the sum $\sum_{i\in
\Z}\mu(\alpha^{-1}(i))\LL^{-i}$ is well-defined in
$\hat{\mathcal{M}}_{k}$. In this case, this sum is by definition
the motivic integral $\int_{A}\LL^{-\alpha}d\mu$.

An important tool in this setting is the change of variables
formula \cite{DLinvent}. Let $X,\,Y$ be varieties over $k$, of
pure dimension $e$, $Y$ smooth, and let $h:Y\rightarrow X$ be a
proper birational morphism. Let $A$ be a measurable subset of
$\mathcal{L}(X)$, and let $\alpha:A\rightarrow \Z$ be a function
such that $\LL^{-\alpha}$ is integrable. Then
\[\int_{A}\LL^{-\alpha}d\mu=\int_{h^{-1}(A)}\LL^{-\alpha\circ h -
ord_{t}h^{*}(\Omega^{e}_{X})}d\mu\ .\] As can be expected, this
formula is often used when $h$ is a resolution of singularities.
It allows one to introduce new invariants of $X$, in terms of a
resolution of singularities, which are independent of the chosen
resolution, since the definition as a motivic integral on
$\mathcal{L}(X)$ is intrinsic. To give an example: using motivic
integration, and the change of variables formula, one can prove
that the topological zeta function associated to a regular
function $f$ is independent of the chosen resolution for $f$
\cite{DL5}.

Let $V\subset X$ be varieties over $k$, $X$ smooth of dimension
$e$, and let $O$ be a closed point of $V$. The embedding of $V$
into the smooth ambient variety $X$ allows us to describe jets on
$V$ by means of arcs on $X$. Let $d$ be a positive integer. We
will compute the local Igusa Poincar\'e series $Q_{geom}(T)$ of
$V$ at $O$ by means of the local motivic Igusa zeta function $Z$,
making use of the change of variables formula. By definition,
\[Z(d)=\int_{\mathcal{L}(X)_{O}}\LL
^{-ord_{t}\mathcal{I}^{d}}d\mu\,,\]
 where
$\mathcal{I}$ is the defining ideal sheaf of $V$ in $X$,  and we
write $\mathcal{L}(X)_{O}$ to denote the arcs on $X$ with origin
in $O$. Recall that, for an arc $\psi$ in $\mathcal{L}(X)_{O}$,
the order $ord_t \mathcal{I}^{d}$ is defined as $\min \{ord_t
f(\psi)\,|\,f\in \mathcal{I}^{d}_{O},\,f(\psi)\neq 0\}$.
 Putting
$T=\LL^{-d}$, we get the classical transformation formula
\[Q_{geom}(T\LL^{-e})=\frac{1 -  \LL^{e}Z(T)}{1-T}\ .\]
If $h:X'\rightarrow X$ is any proper birational morphism, with
$X'$ smooth,
\[Z(d)=\int_{\mathcal{L}(\tilde{X})_{h^{-1}(O)}}\LL
^{-ord_{t}\mathcal{I}^{d}\circ h -ord_t Jac_h}d\mu\,\] where
$Jac_h$ is the Jacobian of $h$. We will take for $h$ an embedded
resolution of $V$ in $X$, because in this case, the behaviour of
$ord_{t}\mathcal{I}^{d}\circ h$ can be made explicit, allowing us
to compute the latter motivic integral.
\section{The motivic Igusa Poincar\'e series}\label{igusa}
Let $V$ be a singular affine toric surface defined by the cone
$\sigma$ generated by $(1,0)$ and $(p,q)$, where $0<p<q$ and $p,q$
are relatively prime. Let $(b_1,\ldots,b_s)$ be the entries
occurring in the Hirzebruch-Jung continued fraction associated to
$q/(q-p)$, and $(c_1,\ldots,c_t)$ the components of the continued
fraction of $q/p$ \cite{F}\cite{Oda}. Let furthermore $\Theta$ be
the union of compact faces of the convex hull of $\sigma\cap
N\setminus 0$, and $\check{\Theta}$ be the union of compact faces
of the convex hull of $\check{\sigma}\cap M\setminus 0$.

The minimal resolution of $V$ is a toric modification induced by a
subdivision of $\sigma$ into simple cones. The vectors occurring
in this subdivision can be listed as follows:
\[v_0=(1,0),\,v_1=(1,1),\ldots,\,v_{j+1}=b_jv_{j}-v_{j-1},\ldots,\,v_{s+1}=b_sv_{s}-v_{s-1}=(p,q)\,.\]
The exceptional divisors $E_j\cong \Pro^{1}$ of this resolution,
$j=1,\ldots,s$, correspond to the newly introduced vectors $v_j$,
and $E_j$ is known to have self-intersection number $-b_j$.

The $c_j$ have a geometric significance of their own: subdividing
$\check{\sigma}$ into simple cones, i.e. taking the minimal set of
generators for the semi-group $\check{\sigma}\cap M$, yields an
embedding of $V$ into affine $(t+2)$-space; the ideal of $V$ is
generated by $x_{i-1}x_{i+1}-x_i^{c_i}$, $i=1,\ldots,t$.

One can derive the $b_i$ from the $c_j$, as will be proved using
the polar polyhedron $\Theta^{0}$ associated to $\Theta$. We will
describe this connection algorithmically. Read the $b_i$ by order
of indexing; a sequence of $j$ 2's induces a dual component $j+3$,
unless this sequence contains $b_1$ or $b_s$; in that case, the
induced number is $j+2$ if only one of both is included, and $j+1$
if the sequence includes both $b_1$ and $b_s$. A value $b_i\neq 2$
induces a dual sequence of $b_i-3$ 2's, unless $i=1$ or $i=s$; in
that case, $b_i-2$ 2's appear, or only $b_i-1$ if $s=1$. If the
successor of $b_i\neq 2$ again differs from 2, these dual 2's must
be followed by a 3. Now move on to $b_{i+1}$ and repeat the
procedure.

\begin{lemma}
Let $V$ be a singular affine toric surface defined by the cone
$\sigma$ generated by $(1,0)$ and $(p,q)$, where $0<p<q$ and $p,q$
are relatively prime. Let $(b_1,\ldots,b_s)$ be the entries
occurring in the Hirzebruch-Jung continued fraction associated to
\\ $q/(q-p)$, and $(c_1,\ldots,c_t)$ the components of the continued
fraction of $q/p$. The algorithm described above computes, with
input $(b_1,...,b_s)$, the output string $(c_1,...,c_t)$.
\end{lemma}

The algorithm works in both ways, that is, also allows one to
deduce the $b_i$ from the $c_j$.

\begin{proof}
Consider the support function $\check{h}:\check{\sigma}\rightarrow
\R^{+}$ of $\Theta$, mapping a vector $m$ in $\check{\sigma}$ to
\[\check{h}(m)=\min \,\{<m,n>\,|\,n\in \Theta \}\ .\] We define
the polar polyhedron $\Theta^{0}$ for $\Theta$ by
\[\Theta^{0}=\{m\in \check{\sigma}\,|\,\check{h}(m)\geq 1\,\}\ .\]
Let $v_{j(\alpha)}$, $\alpha=0,\ldots,l$, be the vertices of
$\Theta$, with $0=j(0)<\ldots<j(l)=s+1$. Let $m_{j(0)}$ be the
vector $(0,1)$ in $\check{\sigma}$, and let $m_{j(l+1)}$ be
$(q,-p)$. We define $m_{j(\alpha)}$, for $\alpha=1,\ldots,l$, to
be the primitive vector in $\check{\sigma}\cap M$ satisfying
\[<m_{j(\alpha)},v_{j(\alpha-1)}>=<m_{j(\alpha)},v_{j(\alpha)}>=1\ .\]
Now $\check{\Theta}$ is the convex hull of $\{m_{j(\alpha)}\,|\,
\alpha=0,\ldots,l+1\,\}+\check{\sigma}$. For $1\leq \alpha\leq
l-1$, the line segment joining $m_{j(\alpha)}$ and
$m_{j(\alpha+1)}$ contains exactly $b_{j(\alpha)}-1$ lattice
points. If we denote by $\{m,m'\}$ the $\Z$-basis for $M$, dual to
$\{v_{j(\alpha)-1},v_{j(\alpha)}\}$, then $m_{j(\alpha)}=m+m'$,
and $m_{j(\alpha+1)}=(b_{j(\alpha)}-1)m+m'$. This information
allows you to derive the algorithm. As an example, let us compute
$c_{1}$, assuming that $j(1)\neq s$. Since the segment joining
$m_{j(1)}=(1,0)$ and $m_{j(2)}$ contains $b_{j(1)}-1$ lattice
points , we find \[(c_1 (b_{j(1)}-2)-b_{j(1)}+3,2-b_{j(1)})\] as
coordinates of $m_{j(2)}$. On the other hand, the basis dual to
$\{v_{j(1)-1},v_{j(1)}\}$ is $\{(j(1),-1),(1-j(1),1)\}$. Comparing
the two expressions for $m_{j(2)}$ yields $c_{1}=j(1)+1$. Now
$j(1)-1$ is the number $j$ of 2's at the beginning of the series
$b_1,\ldots,b_s$. One can show that $c_{b_{j(1)}-1}\neq 2$: thus
the series $c_1,\ldots,c_t$ starts with $b_{j(1)}-2$ 2's, if
$j(1)=1$, and with the number $j+2$ else.
\end{proof}

In order to compute the Igusa Poincar\'e series by means of a
motivic integral and the change of variables formula, we need to
embed the minimal toric resolution for $V$ into an embedded
resolution for $V$ in some smooth ambient space. We will factor
the canonical toric resolution into a sequence of blow-ups of
zero-dimensional orbits, which can be immediately extended to an
embedded resolution for $V$ using the embedding in affine space
mentioned above. Blowing up the unique zero-dimensional orbit $O$
of $V$ corresponds, by \cite{Kempf}, to the toric modification
corresponding to the subdivision $\Sigma$ of $\sigma$ introducing
all  primitive vectors normal to the edges of $\check{\Theta}$.
Using $\Theta^{0}$ to describe $\check{\Theta}$, one can show that
this comes down to inserting $v_1$, $v_{s-1}$, and all $v_i$
determining vertices of $\Theta$, i.e. the $v_i$ for which
$b_i\neq 2$. This is a resolution if and only if all $b_i$ with
$i\neq 1,s$ are different from 2; in the other case, we have to
blow up some more.

The singularities left after blowing up $O$ are all rational
double points of type $A_{c}$. In fact, they are recovered from
the $b_i$ by omitting $b_1$ and $b_s$, and isolating all sequences
of 2's in the remaining $b_i$. Let $c$ be the number of 2's in
such a sequence. This number $c$ can be recovered from the $c_j$:
it is equal to $c_j-3$, with $j$ chosen such that the vertex of
$\check{\Theta}$ corresponding to $x_j$ lies on the two edges
whose normal directions determine the cone in our fan $\Sigma$
corresponding to this sequence of 2's. Moreover, each of the $c_j$
which is bigger than 3 will induce a singularity in this way. The
singularity will be resolved after blowing up the zero-dimensional
orbit corresponding to the associated singular cone (thus
inserting 2 vectors, or 1 if $c=1$) and repeating this procedure
$\lfloor c/2 \rfloor$ times.

This factorization allows us to embed our resolution in ambient
affine space, simply by blowing up the corresponding points in
this space. Let $h:\tilde{X}\rightarrow X=\A_{k}^{t+2}$ be the
proper birational morphism obtained in this way, and let
$\tilde{V}$ be the strict transform of $V$ (thus $\tilde{V}$ is
the canonical resolution surface). The points of $\tilde{V}$ where
there's no transversal intersection with the exceptional locus of
$h$ correspond to adjacent vectors in the simple subdivision of
$\sigma$ which are introduced in one and the same blow-up. In
these points, the intersection multiplicity will be two. This
tells us which jets are shared by $\tilde{V}$ and the exceptional
locus, preventing us from counting them double.

 Let $d$ be a positive integer. We will compute the local Igusa Poincar\'e series $Q_{geom}(T)$ of $V$ at $O$
 by means of the local
 motivic Igusa zeta function $Z$, making use
of the change of variables formula
\[Z(d)=\int_{\mathcal{L}(X)_{O}}\LL ^{-ord_{t}\mathcal{I}^{d}}d\mu
=\int_{\mathcal{L}(\tilde{X})_{h^{-1}(O)}}\LL
^{-ord_{t}\mathcal{I}^{d}\circ h -ord_t Jac_h}d\mu\,\] where
$\mathcal{I}$ is the defining ideal sheaf of $V$ in $X$, and
$Jac_h$ is the Jacobian of $h$. Putting $T=\LL^{-d}$, we get that
\[Q_{geom}(T\LL^{-(t+2)})=\frac{1 -  \LL^{t+2}Z(T)}{1-T}\ .\]
Observe that we can recover $t$ from $Q_{geom}(T)$, since the
coefficient of the $T$-term in the series equals $\LL^{t+2}$.
\section{The computations}\label{comput}
Let $a$ be the number of vectors introduced in $\Sigma$, i.e. the
number of elements in $\{b_2,\ldots,b_{s-1}\}$ differing from 2
augmented by two, and let $b=a-r-1$ be the number of pairs of
adjacent vectors in $\Sigma$, that is, pairs of vectors in
$\Sigma$ with multiplicity $1$.

 We split $\LL^{t+2}Z$ up in
different terms, corresponding to the classical stratification of
the exceptional locus. Let $E$ be the strict transform of the
exceptional divisor that is created by blowing up $O$. The
contribution of arcs in $\tilde{X}$ with origin in $E$, but not in
another exceptional divisor or $\tilde{V}$, is clearly equal to
\[Z_1(d)=([\Pro^{t+1}]
-a[\Pro^{1}]+(a
-1))\frac{(\LL-1)\LL^{-2d-t-2}}{1-\LL^{-2d-t-2}}\,.\] Next, we
consider arcs with origin in the smooth part of $E'=E\cap
\tilde{V}$. In these points, $E$ and $\tilde{V}$ intersect
transversally. Let us denote this set of origins by $E'^{o}$.
Since the order of $\mathcal{I}\circ h$ on $E$ equals 2, in each
point
 the contribution of arcs tangent to neither $E$ nor
$\tilde{V}$ amounts to
\[\alpha:=(\LL^{t+2}-\LL^{t+1}-\LL^{2}+\LL)\LL^{-3d-2t-3}\,.\]
Counting the arcs tangent to $\tilde{V}$ but not to $E$ yields
\[\beta:=(\LL-1)(\LL^{t}-1)\frac{\LL^{-4d-3t-2}}{1-\LL^{-d-t}}\,,\]
while the arcs tangent to $E$ but not to $\tilde{V}$ contribute
\[\gamma:=(\LL-1)(\LL^{t}-1)\frac{\LL^{-5d-3t-4}}{1-\LL^{-2d-t-2}}\,.\]
As for the arcs tangent to both $E$ and $\tilde{V}$, we get the
same computations at the level of 2-jets, and we obtain
\[\sum_{j=1}^{\infty}(\alpha+\beta+\gamma)\LL^{(-3d-2t-2)j}\] which brings the total
contribution of $E'^{o}$ to
\[Z_2(d)=(a[\Pro^{1}]-2(a-1))\frac{(\LL^{t}-1)(\LL-1)\LL^{-3d-2t-2}}{(1-\LL^{-2d-t-2})(1-\LL^{-d-t})}.\]

We also have to cope with exceptional divisors emerging during the
remainder of the resolution process. The situation is as follows:
after blowing up the origin in $X$, some singularities may remain,
situated in the intersection points of exceptional divisors of
$\tilde{V}$. They are described by the numbers $d_k=c_k-3$,
indicating the length of the corresponding sequence of 2's in the
continued fraction series of $q/(q-p)$. It takes $\lceil d_k/2
\rceil$ blow-ups to resolve them; if $d_k$ is odd, we get a chain
of exceptional divisors intersecting $\tilde{V}$ transversally, if
$d_k$ is even we get an intersection point of multiplicity 2 in
the last stage of the resolution process, and we have to blow up
one of the intersection curves to remedy this situation. So for
each $c_k>3$, we get, after simplification, a contribution
\begin{eqnarray*}
Z^{(k)}(d)&=&(\LL-1)\{\sum_{j=2}^{\lceil d_k/2\rceil}\{\,
(\LL^{t+1}-2\LL+1)\frac{\LL^{-N_j'd-\nu'_{j} }}{(1-\LL^{-N_j'
d-\nu'_j })(1-\LL^{-N_{j+1}' d-\nu'_{j+1} })}
\\&&+2\frac{(\LL-1)(\LL^{t}-1)\LL^{-(N_j'  +1)d-\nu'_j  -t}}{(1-\LL^{-d-t})(1-\LL^{-N_j' d-\nu'_j })(1-\LL^{-N'_{j+1} d-\nu'_{j+1} })}
\,\}
\\&&+([\Pro^{t}]-2)(\LL-1)\frac{\LL^{-6d-3t-5}}{(1-\LL^{-2d-t-2 })(1-\LL^{-4d-2t-3 })}
\\&&+2\frac{(\LL-1)(\LL^{t}-1)\LL^{-7d-4t-5}}{(1-\LL^{-d-t})(1-\LL^{-2d-t-2 })(1-\LL^{-4d-2t-3
})}
\\&&+T^{(k)}\,\}\ ,
\end{eqnarray*}
where $N_j'=2j$ and $\nu'_j=j(t+1)+1$. The expression for
$Z^{(k)}-(\LL-1)T^{(k)}$ can be further simplified to $(\LL-1)$
times
\begin{eqnarray*}
\frac{(\LL^{-d-t}+\LL^{-d+1}-2\LL^{-d}+\LL^{t+1}-2\LL+1)(\LL^{-(\lceil
\frac{d_k}{2}\rceil +2)(2d+t+1)-2}+\sum_{j=2}^{\lceil
\frac{d_k}{2}\rceil}\LL^{-j(2d+t+1)-1})}{(1-\LL^{-d-t})(1-\LL^{-2d-t-2})(1-\LL^{-(\lceil
\frac{d_k}{2} \rceil +1 )(2d+t+1)-1})}.
\end{eqnarray*}
The term $T^{(k)}$ depends on the parity of $d_k$. If $d_k$ is
odd, then the exceptional divisor created in the final blow-up
intersects $\tilde{V}$ transversely along a $\Pro^{1}$, so
\begin{eqnarray*}
T^{(k)}&=&([\Pro^{t+1}]-[\Pro^{t}]-[\Pro^{1}]+2)\frac{\LL^{-(d_k+3)
d-((d_k+1)/2+1)(t+1)-1}}{1-\LL^{-(d_k+3) d-((d_k+1)/2+1)(t+1)-1 }}
\\&&+([\Pro^{1}]-2)\frac{(\LL^{t}-1)\LL^{-(d_k+4)d-((d_k+1)/2+1)(t+1) -t-1}}{(1-\LL^{-(d_k+3) d-((d_k+1)/2+1)(t+1)-1
})(1-\LL^{-d-t})}\,.
\end{eqnarray*}
If $d_k$ is even, then the exceptional divisor of the final
blow-up induces two exceptional divisors in $\tilde{V}$,
intersecting in a point where the intersection multiplicity of
$\tilde{V}$ with the global exceptional divisor is 2. We have to
blow up $\tilde{X}$ along one of these divisors to obtain
transversal intersection. Hence, $T^{(k)}$ equals
\begin{eqnarray*}
&&([\Pro^{t+1}]-[\Pro^{t}]-2[\Pro^{1}]+3)\frac{\LL^{-(d_k+2)
d-(d_k/2+1)(t+1)-1}}{1-\LL^{-(d_k+2) d-(d_k/2+1)(t+1)-1 }}
\\&&+(2[\Pro^{1}]-4)\frac{(\LL^{t}-1)\LL^{-(d_k+3)d-(d_k/2+1)(t+1) -t-1}}{(1-\LL^{-(d_k+2)
d-(d_k/2+1)(t+1)-1 })(1-\LL^{-d-t})}
\\&&+([\Pro^{t}]-[\Pro^{t-1}])\frac{\LL^{-(d_k+3)d-(d_k/2+2)(t+1)}}{1-\LL^{-(d_k+3)d-(d_k/2+2)(t+1)}}
\\
&&+\frac{([\Pro^{t-1}]-1)(\LL-1)\LL^{-(2d_k+5)d-(d_k+3)(t+1)-1}}{(1-\LL^{-(d_k+2)d-(d_k/2+1)(t+1)-1})(1-\LL^{-(d_k+3)d-(d_k/2+2)(t+1)})}
\\ &&+
\frac{(\LL-1)(\LL^{t}-1)\LL^{-(2d_k+6)d-(d_k+4)(t+1)-2}}{(1-\LL^{-d-t})(1-\LL^{-(d_k+2)d-(d_k/2+1)(t+1)-1})(1-\LL^{-(d_k+3)d-(d_k/2+2)(t+1)})}\
.
\end{eqnarray*}

To conclude, let us consider arcs with origin in the set of
singular points of $E'$: these are intersection points of
exceptional divisors of $h|_{\tilde{V}}$ of the first generation.
We know that in these points the intersection multiplicity of
$\tilde{V}$ and $E$ is 2; the tangent plane of $\tilde{V}$ will be
contained in $E$. We blow up irreducible components of $E'$ in
order to remedy this situation. Let $F_1$ and $F_2$ be irreducible
components of $E\cap \tilde{V}$, intersecting in a point $x$. We
are interested in the contribution
\[\LL^{t+2}\int_{\mathcal{L}(\tilde{X})_{x}}\LL^{-ord_{t}\mathcal{I}^{d}\circ h - ord_t Jac_h}d\mu\ .\]
Blowing up $\tilde{X}$ along $F_2$ introduces a $\Pro^{t}$-bundle
$F$ over $F_{2}$. By abuse of notation, we denote the strict
transform of the exceptional divisor $E$ again by $E$, and the
strict transform of $\tilde{V}$ by $\tilde{V}$; $E$, $\tilde{V}$
and $F$ intersect transversally, and $E$ meets the fiber of $F$
over $x$ in a $\Pro^{t-1}$. This space $\Pro^{t-1}$ contains the
fiber of the intersection of $\tilde{V}$ and $F$, which is a
point. Applying the change of variables formula yields
\begin{eqnarray*}
Z_3(d)&=&([\Pro^{t}]-[\Pro^{t-1}])(\LL-1)\frac{\LL^{-3d-2t-2}}{1-\LL^{-3d-2t-2}}
\\ &&+([\Pro^{t-1}]-1)(\LL-1)^{2}\frac{\LL^{-5d-3t-4}}{(1-\LL^{-3d-2t-2})(1-\LL^{-2d-t-2})}
\\ &&+(\LL-1)^{2}(\LL^{t}-1)\frac{\LL^{-6d-4t-4}}{(1-\LL^{-d-t})(1-\LL^{-3d-2t-2})(1-\LL^{-2d-t-2})}\
.
\end{eqnarray*}

Let us summarize these results in a more surveyable way. Let $h'$
be the toric modification which is obtained by taking the
canonical resolution, corresponding to the simple subdivision
$\Sigma_{0}$ of $\sigma$, and blowing up a divisor on $\tilde{V}$
through each point with intersection multiplicity 2 which results
from the resolution of the $A_{c}$-singularities; $h'$ does not
include the blow-ups of divisors through singular points of $E'$.
 Define $E_{-1}$ to be the strict transform
of $V$ under $h'$. Let $E_0$ be the strict transform of the
exceptional divisor that is created in the first blow-up, and let
$E_{i,j}$ be the strict transform of the exceptional divisor
induced by the $j$-th blow-up of the singularity corresponding to
the $i$-th sequence of 2's in $b_2,\ldots,b_{s-1}$.

 We let $I$ denote the index set
\[\{-1,0\}
\cup \{(i,j)\,|\,i\in \{1,\ldots,r\},\,j\in\{1,\ldots,\lceil
(d_i+1)/2\rceil \}\}\,.\] Observe that we allow $j$ to range to
$\lceil (d_i+1)/2\rceil$, because of the extra blow-up we
introduced, if necessary, to cope with the point with intersection
multiplicity 2. We stratify $\tilde{X}$ in the usual way: for each
subset $J$ of $I$, we define $E_{J}$ to be $\cap_{\alpha\in
J}E_{\alpha}$, while $E_{J}^{o}$ denotes $E_{J}\setminus
\cup_{\alpha\notin J}E_{\alpha}$.

We attach to each $E_{\alpha}$ a pair of numerical data
$(N_{\alpha},\nu_{\alpha})$ as follows: \[(N_{-1},\nu_{-1})=(1,t),
\,(N_{0},\nu_{0})=(2,t+2), 
(N_{(i,j)},\nu_{(i,j)})=(2(j+1),(j+1)(t+1)+1)\,.\] If $d_i$ is
even and $j=\lceil (d_i+1)/2 \rceil$, we redefine
$(N_{(i,j)},\nu_{(i,j)})$ as \[(d_i+3,(d_i/2+2)(t+1))\,.\]
 Then
\begin{eqnarray*}
Z(d)&=&\LL^{-(t+2)}\sum_{J\subset
I,J\nsubseteq\{-1\}}[E_{J}^{o}]\prod_{\alpha\in J}
\frac{(\LL^{codim\,E_{\alpha}}-1)\LL^{-{N_{\alpha}
d-\nu_{\alpha}}}}{1 -\LL^{-{N_{\alpha} d-\nu_{\alpha}}}}
\\&& \qquad
+b\,\LL^{-(t+2)}\{Z_3(d)-\frac{(\LL^{t}-1)(\LL-1)\LL^{-3d-2t-2}}{(1-\LL^{-d-t})(1-\LL^{-2d-t-2})}\}
\end{eqnarray*}
 This
formula may be considered as a generalization of the formula in
terms of an embedded resolution with normal crossings in the
hypersurface case. The last term corrects for non-transversal
intersection in the singular points of $E'$.

It follows from results in \cite{DLinvent} that this formula holds
already over $\mathcal{M}_{k}$.

\section{Extracting information from the motivic zeta function}\label{poles}
The question presents itself what information is contained in the
Igusa Poincar\'e series of a toric surface singularity, or,
equivalently, in its motivic zeta function.

\begin{theorem}\label{info}
The motivic Igusa Poincar\'e series $Q_{geom}(T)$ determines the
set $\{c_j\}_{j=1}^{t}$.
\end{theorem}

This is the best we can hope for, since the resolution of the
singularity $O$ is, intuitively speaking, independent of the order
of the $c_j$, modulo cutting and pasting. Theorems 1 and 2 imply
that for toric surface singularities, the motivic Igusa Poincar\'e
series contains more information than the geometric and arithmetic
Poincar\'e series. For instance, Corollary 4.9 in \cite{LejReg}
states that $P_{geom}(T)$ is trivial, i.e. equal to
$1/(1-\LL^{2}T)$, if and only if $s=1$. In this case, the set
$c_1,\ldots,c_t$ will consist entirely of 2's. The geometric
series, only considering liftable jets, can't tell you the value
of the multiplicity $t$; the Igusa Poincar\'e series can, since
$t$ appears already in the dimension of the tangent space, i.e.
the space of 1-jets. It follows from \cite{LejReg}, Corollary 4.8,
that this is the only difference between the series $Q_{geom}(T)$
and $P_{geom}(T)$.

\begin{proof}[Proof of Theorem \ref{info}:]
If $t\neq 1$, we can list the candidate poles of the zeta function
$Z(d)$ as follows:
\begin{eqnarray*}&&-t\leq -\frac{2t+2}{3}\leq
-\frac{t+2}{2}<-\frac{(j+1)(t+1)+1}{2(j+1)},\,j\in \{\lceil
d_k/2\rceil\,|\,c_k>3\},\\&&\qquad\qquad
-\frac{(d_k/2+2)(t+1)}{d_k+3},\,d_k\in 2\Z\,.\end{eqnarray*} It is
important for our purposes that
\[\frac{(d_k/2+2)(t+1)}{d_k+3}>\frac{(d_k/2+1)(t+1)+1}{d_k+2}\] if
$t\neq 1$.

The candidate pole $d=-t$ will always be an actual pole of the
motivic zeta function, since in the other case, the denominator
$1-\LL^{2}T$ would not appear in $Q(T)$\,; but if we evaluate
$(1-\LL^{2}T)Q(T)$ at $T=\LL^{-2}$, we get
$[\mathcal{L}_{n}(X)_{O}]\LL^{-2n}$ as the $n$-th partial sum,
i.e. as the evaluation at $T=\LL^{-2}$ of $(1-\LL^{2}T)Q(T)\
\mathrm{mod}\,T^{n+1}$. It follows from \cite{DLinvent}, Theorem
7.1, that
\[\lim_{n\to
\infty}[j^{n}\mathcal{L}(X)_O]\LL^{-2n}=\LL^{2}\mu(\mathcal{L}(X)_{O})\,,\]
which is non-zero; so the series $[\mathcal{L}_{n}(X)_O]\LL^{-2n}$
diverges, or has a nonzero limit.
Hence, we recover $t$ by looking at the smallest pole of $Z$,
which is $-t$ or non-integer - in the latter case, $t$ must be
equal to one.

Let us investigate what happens when we specialize to the
topological zeta function, as is described in \cite{DL5}. This
yields
\begin{eqnarray*}
Z_{top}(d)&=&\sum_{J\subset I}\chi[E_{J}^{o}]\prod_{\alpha\in J}
\frac{1}{N_{\alpha}d+\nu_{\alpha}}
\\&& \qquad
+b\,\{Z_{3,\chi}(d)-\frac{1}{(d+t)(2d+t+2)}\}\,,
\end{eqnarray*}
where $\chi:\Z[\LL]\rightarrow \Z$ is the topological Euler
characteristic; this simply means that we write all coefficients
in terms of $\LL$, and map $\LL$ to 1. The function $Z_{top}(d)$
is well-defined, since it is a specialization of the motivic zeta
function, which is defined intrinsically. Poles of $Z_{top}$ will
correspond to poles of $Z$, since the Grothendieck bracket is a
finer invariant than the Euler characteristic. Working with
$Z_{top}$ instead of $Z$ obviously simplifies the computations,
but when the Euler characteristic is too coarse to detect certain
poles, or to give useful information about their residues, one is
obliged to turn back to $Z$.

The case $t=1$ being trivial, we might as well assume that $t>1$.
First suppose $t\neq 2$. The residue of the candidate pole
$d=-(2t+2)/3$ is equal to
\[-\frac{b}{3(t-2)^{2}}\{2t^{2}-5t+11\}\,\]
 which enables us to recover
the value of $b$.
 If $t$
happens to be 2, we can still recover $b$ by looking at the
residue of the pole $d=-2$, which will have multiplicity 3.

The largest candidate pole of $Z_{top}$ is the one induced by
$(N_{i,j},\nu_{i,j})$, with $d_i$ maximal among the $d_k$,
 and $j=\lceil d_i/2 \rceil$. Its residue depends on the numbers $\delta$, $\epsilon$ of occurrences of
 $2j$, resp. $2j-1$,
among the $d_k$.
 At any rate, it is strictly positive if $r\neq 0$, since it concerns the residue of the largest candidate
 pole, and all relevant $\chi[E_{J}^{o}]$ are positive.
Hence, by looking at the largest pole, and its residue with
respect to $Z_{top}$, we can determine $\lceil d_i/2 \rceil$, and
we get a linear relation on $\delta$ en $\epsilon$. An additional,
independent linear relation is obtained by studying the evaluation
of
\begin{eqnarray*}&&Z(d)(1-\LL^{-2(j+1)d-(j+1)(t+1)-1})(1-\LL^{-d-t})(1-\LL^{-2jd-j(t+1)-1})
\\&&\qquad\qquad\qquad\qquad\qquad\qquad\qquad
\qquad\qquad\qquad .(1-\LL^{-(2j+3)d-(j+2)(t+1)})\end{eqnarray*}
at $d=-((j+1)(t+1)+1)/2(j+1)$. Since we are considering the
largest candidate pole, one sees without further calculations that
the coefficient of $\LL^{t+2}$ is equal to $\delta+\epsilon$.
Working backwards to $-(2t+3)/4$, we can determine all $d_k$.
\end{proof}

\noindent \textit{Remark:} Since we don't know wether
$\mathcal{M}_k$ is a domain, we should explain what we mean by a
pole of a rational function over $\mathcal{M}_k$. An exact
definition is given in \cite{RoVe}.

\bibliographystyle
{hplain}
\bibliography{wanbib,wanbib2}

\end{document}